\documentclass[11pt]{article}
\usepackage{enumerate}
\usepackage{amssymb,a4wide,latexsym,makeidx,epsfig,fleqn}
\usepackage{amsthm}
\usepackage{amsmath}
\usepackage{enumerate}
\newtheorem{theorem}{Theorem}[section]

\newtheorem{lemma}[theorem]{Lemma}

\newtheorem{corollary}[theorem]{Corollary}

\begin{document}
\textwidth 150mm \textheight 225mm
\title{Signless Laplacian spectral conditions for Hamilton-connected graphs with large minimum degree
\thanks{ Supported by
the National Natural Science Foundation of China (No. 11171273)}}
\author{{Qiannan Zhou$^1$, Ligong Wang$^{2,}$\footnote{Corresponding author.} and Yong Lu$^3$}\\
{\small Department of Applied Mathematics, School of Science, Northwestern
Polytechnical University,}\\ {\small  Xi'an, Shaanxi 710072,
People's Republic
of China.}\\
{\small $^1$ E-mail: qnzhoumath@163.com}\\
{\small $^2$ E-mail: lgwangmath@163.com}\\
{\small $^3$ E-mail: luyong.gougou@163.com}\\}
\date{}
\maketitle
\begin{center}
\begin{minipage}{120mm}
\vskip 0.3cm
\begin{center}
{\small {\bf Abstract}}
\end{center}
{\small In this paper, we present a spectral sufficient condition for a graph to be Hamilton-connected in terms
of signless Laplacian spectral radius with large minimum degree.

 }
\end{minipage}
\end{center}

\section{Introduction }

Throughout this paper, we only consider simple graphs without loops and multiple edges. Let $G=(V(G),E(G))$ be a graph with vertex set $V(G)=\{v_{1},v_{2},\ldots,v_{n}\}$ and edge set $E(G)$. Denote by $e(G)=|E(G)|$ the edge number of $G$. Let $V_{1}$ and $V_{2}$ be disjoint subsets of $V(G)$. We denote by $E(V_{1},V_{2})$ the set of edges each of which has one vertex in $V_{1}$ and the other vertex in $V_{2}$ and let $e(V_{1},V_{2})=|E(V_{1},V_{2})|$.
The degree of $v$ is denoted by $d_{G}(v)=|N_{G}(v)|$, where $N_{G}(v)$ is the set of vertices adjacent to $v$ in $G$. Moreover, $N_{G}[v]=N_{G}(v)\cup \{v\}$. Denote by $\delta(G)$ the minimum degree of $G$ and $\omega(G)$ the clique number of $G$. If the graph $G$ is clear under the context, we will drop the subscript $G$. We use
$G+H$ and $G\vee H$ to denote the disjoint union and the join of $G$ and $H$, respectively. The union of $k$ disjoint
copies of the same graph $G$ is denoted by $kG$. For an edge subset $E_{0}$ of a graph $G$, $G-E_{0}$ denotes the graph obtained from $G$ by deleting the edges in $E_{0}$. For an induced subgraph $H$ of $G$, denote by $G-H$, the subgraph obtained from $G$ by deleting all vertices of $H$ and all incident edges.  For terminology and notation not
defined but used, we refer the reader to \cite{BondyM2008}.

Let $A(G)$ and $D(G)$ be the adjacency matrix and the degree diagonal matrix of $G$, respectively. The matrix
$Q(G)=D(G)+A(G)$ is called to be the signless Laplacian matrix of $G$. The largest eigenvalue $q(G)$ of $Q(G)$ is
called to be the signless Laplacian spectral radius of $G$.

 A graph is called to be Hamiltonian (traceable) if it contains a Hamilton cycle (path), which is a cycle (path) containing all vertices of $G$. A graph is called to be Hamilton-connected if every two
 vertices of $G$ are connected by a Hamilton path.

It is well-known that the problem of determining the Hamiltonicity of graphs is NP-hard. So scholars focus on
finding sufficient conditions for graphs to be Hamiltonian, traceable or Hamilton-connected. In 2010, Fiedler and Nikiforov \cite{Fiedler2010} gave sufficient conditions on spectral radius for the existence of Hamilton cycles.
Motivated by this, there are many other spectral conditions for the Hamiltonicity of graphs, one may refer to
\cite{Benediktovich2015,Liu2015,Lu2012,NingGe2015,fy,YuYeCai2014,ZhouBo2010,ZhouQN22017,ZhouQN12017}.
Recently, by imposing the minimum degree of a graph as a new parameter, Li and Ning \cite{LiNing2016LMA,LiNing2017LAA} extended some the results in \cite{Fiedler2010,Liu2015,NingGe2015}. Now, their
results were improved by Nikiforov \cite{Nikiforov2016}, Chen et al. \cite{ChenHou2017}, Ge et al. \cite{GeNing} and Li et al. \cite{LiLiuPeng}, in some sense. In this paper, we will establish a signless Laplacian analogue of a result due to Nikiforov \cite{Nikiforov2016} for Hamilton-connected graphs with large minimum degree, which generalizes the result in \cite{ZhouQN12017}.

Firstly, for $n\geq 5$ and $1\leq k\leq n/2$, we define:
$$S_{n}^{k}=K_{k}\vee (K_{n-(2k-1)}+(k-1)K_{1})~\mbox{and}~T_{n}^{k}=K_{2}\vee (K_{n-(k+1)}+K_{k-1}).$$
The graph $S_{n}^{k}$ is obtained from $K_{n-k+1}$ and $(k-1)K_{1}$ by connecting all vertices of $(k-1)K_{1}$ to
all vertices of a $k$-subset of $K_{n-k+1}$. The graph $T_{n}^{k}$ is obtained from $K_{n-k+1}$ and $K_{k+1}$ by
identifying two vertices.

For the graph $S_{n}^{k}$, let $X=\{v\in V(S_{n}^{k}): d(v)=k\}$, $Y=\{v\in V(S_{n}^{k}): d(v)=n-1\}$, $Z=\{v\in V(S_{n}^{k}): d(v)=n-k\}$. Let $E_{0}$ be the subset of $E(S_{n}^{k})$ containing the edges whose both endpoints
are from $Y\cup Z$. For the graph $T_{n}^{k}$, let $X=\{v\in V(T_{n}^{k}): d(v)=k\}$, $Y=\{v\in V(T_{n}^{k}): d(v)=n-1\}$, $Z=\{v\in V(T_{n}^{k}): d(v)=n-k\}$. Let $E_{0}$ be the subset of $E(T_{n}^{k})$ containing the edges whose both endpoints are from $Y\cup Z$. For any real number $x$, let $\lfloor x\rfloor$ denote the greatest integer that is less than or equal to $x$. The integer $\lfloor x\rfloor$ is called the floor of $x$. Then we define:
$$\mathcal{S}_{k}^{(1)}(n)=\{G\in S_{n}^{k}: G=S_{n}^{k}-E',~\mbox{where}~E'\subseteq E_{0}~\mbox{with}~|E'|\leq \lfloor \frac{k(k-1)}{4}\rfloor\}.$$

$$\mathcal{T}_{k}^{(1)}(n)=\{G\in T_{n}^{k}: G=T_{n}^{k}-E',~\mbox{where}~E'\subseteq E_{0}~\mbox{with}~|E'|\leq \lfloor \frac{k-1}{2}\rfloor\}.$$

We have our main theorem:
\noindent\begin{theorem}\label{th:6c0} Let $G$ be a graph of order $n\geq k^{4}+5k^{3}+2k^{2}+8k+12$ with minimum degree $\delta(G)\geq k$, where $k\geq 2$. If
$$q(G)\geq 2n-2k,$$
then $G$ is Hamilton-connected unless $G\in \mathcal{S}_{k}^{(1)}(n) \cup \mathcal{T}_{k}^{(1)}(n)$.
\end{theorem}

Then we give the definitions of another two families of subgraphs of $S_{n}^{k}$ and $T_{n}^{k}$, respectively.

$$\mathcal{S}_{k}^{(2)}(n)=\{G\in S_{n}^{k}: G=S_{n}^{k}-E',~\mbox{where}~E'\subseteq E_{0}~\mbox{with}~|E'|= \lfloor \frac{k(k-1)}{4}\rfloor+1\}.$$

$$\mathcal{T}_{k}^{(2)}(n)=\{G\in T_{n}^{k}: G=T_{n}^{k}-E',~\mbox{where}~E'\subseteq E_{0}~\mbox{with}~|E'|= \lfloor \frac{k-1}{2}\rfloor+1\}.$$

The rest paper is organized as follows. In Section~2, we will give some useful techniques and lemmas. In Section~3,
we present the main theorems, lemmas, and the proofs of them. In Section~4, we give the specific proof of the inequality \eqref{eq:2} which appears in Section~3.

\section{Preliminaries}

By Rayleigh's principle, we have
$$q(G)=\max\limits_{\textbf{x}}\frac{\langle Q(G)\textbf{x},\textbf{x}\rangle}{\langle \textbf{x},\textbf{x}\rangle},$$
where $q(G)$ is the largest eigenvalue of $Q(G)$. If $\textbf{f}$ is the eigenvector corresponding to $q(G)$, then
we get $\textbf{f}_{v}> 0$ for each $v\in V(G)$ by the famous Perron-Frobenius theorem \cite{GodsilRoyle2001}. It is
easy to get
\begin{equation}\label{eq:0}
(q(G)-d(v))\textbf{f}_{v}=\sum_{u\sim v}\textbf{f}_{u}.
\end{equation}

For an integer $k$, the $k$-closure of a graph $G$ is the graph obtained from $G$ successively joining pairs of
nonadjacent vertices whose degree sum is at least $k$ until no such pair remains. We denote it by $cl_{k}(G)$.
Note that $d_{cl_{k}(G)}(u)+d_{cl_{k}(G)}(v)\leq k-1$ for any pair of nonadjacent vertices $u,v\in cl_{k}(G)$ and
$G\subseteq cl_{n+1}(G)$.

\noindent\begin{lemma}\label{le:6c0}(\cite{Ore1960}) If $G$ is a 2-connected graph of order $n$ and $d(u)+d(v)\geq n+1$ for any two distinct nonadjacent vertices $u$ and $v$, then $G$ is Hamilton-connected.
\end{lemma}

\noindent\begin{lemma}\label{le:6c1}(\cite{BondyChvatal1976}) A graph $G$ of order $n$ is Hamilton-connected if and
only if $cl_{n+1}(G)$ is Hamilton-connected.
\end{lemma}

Then let us recall Kelmas transformation \cite{Kelmans1981}. Given a graph $G$ and two specified vertices
$u,v\in V(G)$ construct a new graph $G^{*}=G-\{vx: x\in N_{G}(v)\setminus N_{G}[u]\}+\{ux: x\in N_{G}(v)\setminus N_{G}[u]\}$. The vertices $u$ and $v$ are adjacent in $G^{*}$ if and only if they are adjacent in $G$. Then we have
the following lemma.

\noindent\begin{lemma}\label{le:6c2}(\cite{LiNing2016LMA}) Let $G$ be a graph and $G^{*}$ be
a graph obtained from $G$ by a Kelmans transformation. Then $q(G)\leq q(G^{*})$.
\end{lemma}


The following theorem gives an upper bound for $q(G)$.

\noindent\begin{theorem}\label{th:6cq}(\cite{LFGU,fy}) Let $G$ be a connected graph on $n$ vertices and $m$ edges.
Then $q(G)\leq \frac{2m}{n-1}+n-2$.
\end{theorem}

\section{Main results}

\noindent\begin{theorem}\label{th:6c1} Let $G$ be a graph of order $n\geq 11k$ with minimum degree $\delta(G)\geq k$, where $k\geq 2$. If
$$e(G)> \dbinom{n-k}{2}+k(k+1),$$
then $G$ is Hamilton-connected unless $G\subseteq S_{n}^{k}$ or $T_{n}^{k}$.
\end{theorem}

\noindent {\bf Proof.} Let $G'=cl_{n+1}(G)$. Obviously, we have $\delta(G')\geq \delta(G)$ and $e(G')\geq e(G)$.
If $G'$ is Hamilton-connected, then so is $G$ by Lemma \ref{le:6c1}. Now, we
suppose that $G'$ is not Hamilton-connected. Then we have the following claim.

\noindent\textbf{Claim.} $\omega(G')= n-k+1$.

\begin{proof} Let $K$ be the subset of $V(G')$ containing all vertices which have degree at least $(n+1)/2$. By the
definition of $cl_{n+1}(G)$, any two vertices in $K$ are adjacent in $G'$. Let $C$ be a maximum clique of $G'$
containing all vertices in $K$ and $|C|=t$. Let $H=G'-C$.

Firstly suppose $1\leq t\leq \frac{2}{3}n+\frac{1}{2}$. For every $v\in V(H)$, we have $d_{C}(v)\leq t-1$ and
$d_{G'}(v)\leq \frac{n}{2}$. Note that
$$e(H)+e(V(H),V(C))=\frac{\sum_{v\in V(H)}d_{G'}(v)+\sum_{v\in V(H)}d_{C}(v)}{2}.$$
Hence,
\begin{align*}
e(G')&=e(G'[C])+e(H)+e(V(H),V(C))\\
 &\leq\dbinom{t}{2}+\frac{(\frac{n}{2}+t-1)(n-t)}{2}\\
 &= \frac{n}{4}t+\frac{n(n-2)}{4}\\
 &\leq\frac{n}{4}(\frac{2}{3}n+\frac{1}{2})+\frac{n(n-2)}{4}\\
 &= \frac{5}{12}n^{2}-\frac{3}{8}n\\
 &< \dbinom{n-k}{2}+k(k+1)\\
 &< e(G'),
\end{align*}
which is a contradiction.

Then we suppose that $\frac{2}{3}n+\frac{1}{2}< t\leq n-k$. For every $v\in V(H)$, we have $d_{G'}(v)\leq n-t+1$,
since otherwise $v$ will be adjacent to every vertices in $C$. Note that
$$e(H)+e(V(H),V(C))\leq \sum_{v\in V(H)}d_{G'}(v).$$
Thus,
\begin{align*}
e(G')&=e(G'[C])+e(H)+e(V(H),V(C))\\
 &\leq\dbinom{t}{2}+(n-t+1)(n-t)\\
 &= \frac{3}{2}t^{2}-(2n+\frac{3}{2})t+n^{2}+n\\
 &\leq\dbinom{n-k}{2}+k(k+1)\\
 &< e(G'),
\end{align*}
which is a contradiction.

So we conclude that $t\geq n-k+1$ and $\omega(G')\geq n-k+1$.

Now, suppose that $\omega(G')\geq n-k+2$. Let $C'$ be a maximum clique in $G'$ and $H'=G'-C'$. Since $G'$ is not a
clique, $V(H')\neq \emptyset$. Let $v$ be a vertex in
$V(H')$. Then $v$ is adjacent to every vertex of $C'$ in $G'$ since $d_{G'}(u)\geq n-k+1$ for every $u\in C'$ and
$d_{G'}(v)\geq \delta(G')\geq k$, which contradicts to that $C'$ is a maximum clique in $G'$. So $\omega(G')=
n-k+1$, as desired.
\end{proof}

Note that every vertex in $C'$ has degree at least $n-k$. Let $F=\{u_{1},u_{2},\ldots,u_{s}\}$ be the subset of
$V(C')$ containing the vertices that have degree at least $n-k+1$. In other words, every vertex of $F$ has at least
one neighbor in $H'$. Then we claim that every vertex in $H'$ has degree exactly $k$, which comes from
$G'=cl_{n+1}(G)$ and $\delta(G')\geq k$. Hence every vertex in $H'$ is adjacent to every vertex in $F$. Moreover,
since $|V(H')|=k-1$, we can see that $2\leq s\leq k$.

If $s=2$, then $H'$ is a clique and $G'=T_{n}^{k}$. If $s=k$, then $H'$ is an independent set and $G=S_{n}^{k}$. If $3\leq s\leq k-1$,
we shall show that $G'$ is Hamilton-connected.

Now, we can obviously see that $H'$ is $(k-s)$-regular, $|V(H')|=k-1$ and $G'=K_{s}\vee (K_{n-k-s+1}+H')$, in which
$F=V(K_{s})$ and $C'=K_{s}\vee (K_{n-k-s+1})$.
Let $W=K_{s}\vee H'$. Then $d_{W}(w_{1})+d_{W}(w_{2})=2k\geq s+1+k> |W|+1$ for any two nonadjacent vertices $w_{1}$
and $w_{2}$ in $W$. Hence by Lemma \ref{le:6c0}, $W$ is Hamilton-connected. Then it is obvious that for any two distinct vertices
$u\in V(K_{n-k-s+1})$ and $v\in V(G')$, there is a Hamilton path connecting $u$ and $v$ in $G'$.

For any two distinct vertices $u, v\in F$, we show there is a Hamilton path connecting them in $G'$. Since
every vertex in $H'$ has degree $k-s$ in $H'$, by Dirac's theorem \cite{Dirac1952}, $H'$ has a path
$P=u_{1}u_{2}\cdots u_{t}$ of order $t\geq k-s+1$. Let $V(H')\setminus V(P)=\{v_{1},v_{2},\ldots,v_{k-t-1}\}$. Denote by
$F\setminus \{u,v\}=\{a_{1},a_{2},\ldots,a_{s-2}\}$ and $V(K_{n-k-s+1})=\{b_{1},b_{2},\ldots,b_{n-k-s+1}\}$. If $P$
is a Hamilton path in $H'$, then $uu_{1}\cdots u_{t}a_{1}\cdots a_{s-2}b_{1}\cdots b_{n-k-s+1}v$ is a Hamilton path
connecting $u$ and $v$ in $G'$. If $P$ is not a Hamilton path in $H'$, since $t\geq k-s+1$, then $uu_{1}\cdots
u_{t}a_{1}v_{1}a_{2}v_{2}\cdots a_{k-t-1}v_{k-t-1}a_{k-t}\cdots a_{s-2}b_{1}\cdots b_{n-k-s+1}v$ is a Hamilton path
connecting $u$ and $v$ in $G'$.

For any two distinct vertices $u\in F$, $v\in V(H')$, we show there is a Hamilton path connecting them in $G'$.
Let $F\setminus\{u\}=\{a_{1},a_{2},\ldots,a_{s-1}\}$. If $v\in V(P)$, we use $u^{+}$ and $u^{-}$ to denote the
successor and predecessor of $u$ in $P$, then $uu_{1}\ldots u^{-}a_{s-1}u^{+}\ldots u_{t}a_{1}v_{1}a_{2}v_{2}\cdots
a_{k-t-1}v_{k-t-1}a_{k-t}\\ \cdots a_{s-2}b_{1}\cdots b_{n-k-s+1}v$ is a Hamilton path connecting $u$ and $v$ in $G'$.
If $v\in V(H')\setminus V(P)$, without loss of generality, we assume $v=v_{1}$, then $uu_{1}\ldots u_{t}a_{2}v_{2}\ldots
a_{k-t-1}v_{k-t-1}a_{k-t}\ldots a_{s-1}\\b_{1}\ldots b_{n-k-s+1}a_{1}v$ is a Hamilton path connecting $u$ and $v$ in
$G'$.

For any two distinct vertices $u, v\in V(H')$, we can show there is a Hamilton path connecting them in $G'$ by a
similar method and we omit it.

Hence if $3\leq s\leq k-1$, then $G'$ is Hamilton-connected, which is a contradiction to our assumption.
So $G\subseteq S_{n}^{k}$ or $T_{n}^{k}$.

The proof is complete.
\hfill$\blacksquare$

Then we give the following two lemmas.

\noindent\begin{lemma}\label{th:6c2}  
For each $G\in \mathcal{S}_{k}^{(1)}(n)\cup \mathcal{T}_{k}^{(1)}(n)$, we have $q(G)\geq 2n-2k$.
\end{lemma}

\noindent {\bf Proof.} For $G\in \mathcal{S}_{k}^{(1)}(n)$, we define a column vector $\textbf{c}$ such that
$\textbf{c}_{u}=1$ for all $u\in Y\cup Z$ and $\textbf{c}_{v}=0$ for all $v\in X$. Note that $\textbf{c}$
is the eigenvector corresponding to the eigenvalue $2n-2k$ of $(k-1)K_{1}+ K_{n-k+1}$. Then we have
$$\langle Q(G)\textbf{c},\textbf{c}\rangle-\langle Q((k-1)K_{1}+ K_{n-k+1})\textbf{c},\textbf{c}\rangle=k(k-1)-4|E'|\geq 0.$$
By Rayleigh's principle, we get $q(G)\geq 2n-2k$.

Similarly, for $G\in \mathcal{T}_{k}^{(1)}(n)$, we have
$$\langle Q(G)\textbf{c},\textbf{c}\rangle-\langle Q((k-1)K_{1}+ K_{n-k+1})\textbf{c},\textbf{c}\rangle=2(k-1)-4|E'|\geq 0.$$
By Rayleigh's principle, we can also get $q(G)\geq 2n-2k$.

The proof is complete.
\hfill$\blacksquare$

\noindent\begin{lemma}\label{th:6c3}  Let $G$ be a graph of order $n\geq k^{4}+5k^{3}+2k^{2}+8k+12$ with minimum degree $\delta(G)\geq k$, where $k\geq 2$. For each $G\in \mathcal{S}_{k}^{(2)}(n)\cup \mathcal{T}_{k}^{(2)}(n)$, we have $q(G)<
2n-2k$.
\end{lemma}

\noindent {\bf Proof.} Let $G\in \mathcal{S}_{k}^{(2)}(n)\cup \mathcal{T}_{k}^{(2)}(n)$ and $\textbf{c}$ be
the vector defined in the proof in Lemma \ref{th:6c2}. Then we have the following claim.

\noindent\textbf{Claim~1.} $q(G)> 2n-2k-1$.

\begin{proof} If $G\in \mathcal{S}_{k}^{(2)}(n)$, then we have
$$\langle Q(G)\textbf{c},\textbf{c}\rangle-\langle Q((k-1)K_{1}+ K_{n-k+1})\textbf{c},\textbf{c}\rangle=k(k-1)-4|E'|\geq -4,$$
which implies that $q(G)\geq 2n-2k-\frac{4}{\|\textbf{c}\|^{2}}> 2n-2k-1$.

If $G\in \mathcal{T}_{k}^{(2)}(n)$, then we have
$$\langle Q(G)\textbf{c},\textbf{c}\rangle-\langle Q((k-1)K_{1}+ K_{n-k+1})\textbf{c},\textbf{c}\rangle=2(k-1)-4|E'|\geq -4,$$
which implies that $q(G)\geq 2n-2k-\frac{4}{\|\textbf{c}\|^{2}}> 2n-2k-1$.

Claim~1 is proved.
\end{proof}

In the following proof, we will assume that $G\in \mathcal{S}_{k}^{(2)}(n)$. Since the proofs for the case $G\in \mathcal{T}_{k}^{(2)}(n)$ are similar to those in the case $G\in \mathcal{S}_{k}^{(2)}(n)$, we only give the sketch
in the last.

Let $G\in \mathcal{S}_{k}^{(2)}(n)$ have the maximum signless Laplacian spectral radius. Let $\textbf{f}$ be the
eigenvector corresponding to $q(G)$. Furthermore, we assume $\max_{v\in V(G)}\textbf{f}_{v}=1$.

Since $E'$ is the edge set in which both endpoints are from $Y\cup Z$, we define two subsets of $Y$ and $Z$, respectively. They are as follows:
$$Y_{1}=\{y\in Y: d(y)=n-1\}~\mbox{and}~Y_{2}=\{y\in Y: d(y)\leq n-2\};$$
$$Z_{1}=\{y\in Y: d(y)=n-k\}~\mbox{and}~Z_{2}=\{y\in Y: d(y)\leq n-k-1\}.$$
Since $n\geq k^{4}+5k^{3}+2k^{2}+8k+12$, $Z_{1}\neq \emptyset$. Then we have the following claims.

\noindent\textbf{Claim~2.} $\textbf{f}_{x}\leq \frac{k}{q(G)-k}$ for each $x\in X$.

\begin{proof} By \eqref{eq:0}, we have
$$(q(G)-d(x))\textbf{f}_{x}=\sum_{y\in Y}\textbf{f}_{y}.$$
Since $d(x)=k$, Claim~2 holds.
\end{proof}

\noindent\textbf{Claim~3.}
\begin{enumerate}[(1)]
  \item If $Y_{1}\neq \emptyset$, $Y_{2}\neq \emptyset$, then for any $u\in Y_{1}$, $v\in Y_{2}\cup Z_{1}$, we have $\textbf{f}_{u}> \textbf{f}_{v}$.
  \item If $Z_{2}\neq \emptyset$, $Y_{2}\neq \emptyset$, then for any $u\in Z_{1}$, $v\in Z_{2}\cup Y_{2}$, we have $\textbf{f}_{u}> \textbf{f}_{v}$.
\end{enumerate}

\begin{proof} For any $u, v\in V(G)$, combining with \eqref{eq:0}, we have
\begin{align}\label{eq:1}
(q(G)-d(u))(\textbf{f}_{u}-\textbf{f}_{v})&=\nonumber (q(G)-d(u))\textbf{f}_{u}-(q(G)-d(v))\textbf{f}_{v}+(d(u)-d(v))\textbf{f}_{v}\\
 &\nonumber=(d(u)-d(v))\textbf{f}_{v}+\sum_{s\in N(u)}\textbf{f}_{s}-\sum_{t\in N(v)}\textbf{f}_{t}\\
 &=(d(u)-d(v))\textbf{f}_{v}+\sum_{s\in N(u)\setminus N(v)}\textbf{f}_{s}-\sum_{t\in N(v)\setminus N(u)}\textbf{f}_{t}.
\end{align}

For any $u\in Y_{1}$ and $v\in Y_{2}\cup Z_{1}$; or $u\in Z_{1}$ and $v\in Z_{2}$, then $uv\in E(G)$ and \eqref{eq:1} is equivalent to the following equation:

\begin{equation}\label{eq:11}
(q(G)-d(u)+1)(\textbf{f}_{u}-\textbf{f}_{v})=(d(u)-d(v))\textbf{f}_{v}+\sum_{s\in N(u)\setminus N[v]}\textbf{f}_{s}-\sum_{t\in N(v)\setminus N[u]}\textbf{f}_{t}.
\end{equation}

Since $N(v)\setminus N[u]=\emptyset$, by \eqref{eq:11}, then
$$(q(G)-d(u)+1)(\textbf{f}_{u}-\textbf{f}_{v})=(d(u)-d(v))\textbf{f}_{v}+\sum_{s\in N(u)\setminus N[v]}\textbf{f}_{s}> 0,$$
so $\textbf{f}_{u}> \textbf{f}_{v}$.

Since $G$ has the maximum
signless Laplacian spectral radius in $\mathcal{S}_{k}^{(2)}(n)$, we claim that $G[Y]$ contains the largest number
of edges in $\mathcal{S}_{k}^{(2)}(n)$. If not, since $|E'|=\lfloor \frac{k(k-1)}{4}\rfloor+1$ and $n\geq
k^{4}+5k^{3}+2k^{2}+8k+12$, there always exist $v,w\in Y_{2}$, $u\in Z_{2}$ such that $vw\notin E(G)$ and $uw\in
E(G)$. We construct $G^{*}$ by adding edges $\{vw_{i}: w_{i}\in N_{G}(u)\setminus N_{G}[v]\}$ and deleting edges
$\{uw_{i}: w_{i}\in N_{G}(u)\setminus N_{G}[v]\}$. Hence, $G^{*}[Y]$ has more edges
than $G[Y]$. It is easy to see that the above transformation is the Kelmans transformation. By Lemma \ref{le:6c2},
we have $q(G^{*})\geq q(G)$, which contradicts to choice of $G$.

For any $u\in Z_{1}$ and $v\in Y_{2}$, we shall prove $\textbf{f}_{u}> \textbf{f}_{v}$. Suppose there exist $u\in Z_{1}$ and $v\in Y_{2}$ such that $\textbf{f}_{u}\leq \textbf{f}_{v}$. Let $w\in Y_{2}$ be a vertex not
adjacent to $v$. We construct a new graph $G_{0}$ by adding one edge $vw$ and deleting the edge $uw$. Note that
$$\langle Q(G_{0})\textbf{f},\textbf{f})\rangle-\langle Q(G)\textbf{f},\textbf{f})\rangle=\textbf{f}_{v}^{2}+2\textbf{f}_{v}\textbf{f}_{w}-\textbf{f}_{u}^{2}-2\textbf{f}_{u}
\textbf{f}_{w}\geq 0.$$
We get $q(G_{0})\geq q(G)$ and $G_{0}[Y]$ has more edges than $G[Y]$, a contradiction.

Claim~3 is proved.
\end{proof}

\noindent\textbf{Claim~4.} $\max\limits_{v\in Y\cup Z}\textbf{f}_{v}-\min\limits_{u\in Y\cup Z}\textbf{f}_{u}\leq \frac{k^{2}+6k+6}{2(q(G)-n+1)}.$

\begin{proof} We discuss the following two cases.

\noindent\textbf{Case~1.} $Y_{1}=\emptyset$.

In this case, by Claim~3, $\max_{v\in Y\cup Z}\textbf{f}_{v}$ is attained by some vertex from $Z_{1}$, say, $z$.
By the definition of $Z_{1}$, $z$ is adjacent to other vertices $Y\cup Z$.

If $\min_{v\in Y\cup Z}\textbf{f}_{v}$ is attained by some vertex from $Z_{2}$, say, $w$. Obviously, we have
$zw\in E(G)$, $N(w)\setminus N[z]=\emptyset$, $d(z)-d(w)\leq \lfloor\frac{k(k-1)}{4}\rfloor+1$ and $|N(z)\setminus N[w]|\leq \lfloor\frac{k(k-1)}{4}\rfloor+1$. Hence, by \eqref{eq:11}, we have
\begin{align*}
(q(G)-d(z)+1)(\textbf{f}_{z}-\textbf{f}_{w})&=(d(z)-d(w))\textbf{f}_{w}+\sum_{s\in N(z)\setminus N[w]}\textbf{f}_{s}\\
 &\leq \frac{k(k-1)}{2}+2.
\end{align*}
Since $d(z)=n-k$, we obtain
$$\textbf{f}_{z}-\textbf{f}_{w}\leq \frac{k(k-1)+4}{2(q(G)-n+k+1)}=\frac{k^{2}-k+4}{2(q(G)-n+k+1)}< \frac{k^{2}+6k+6}{2(q(G)-n+1)}.$$

If $\min_{v\in Y\cup Z}\textbf{f}_{v}$ is attained by some vertex from $Y_{2}$, say, $w$. Obviously, we have
$zw\in E(G)$, $N(w)\setminus N[z]=X$ and $|N(z)\setminus N[w]|\leq \lfloor\frac{k(k-1)}{4}\rfloor+1$. Furthermore, we observe
$$d(z)-d(w)\leq ((n-k)-(n-1-\lfloor\frac{k(k-1)}{4}\rfloor-1))=\lfloor\frac{k(k-1)}{4}\rfloor-k+2.$$
Hence by \eqref{eq:11}, we have
\begin{align*}
(q(G)-d(z)+1)(\textbf{f}_{z}-\textbf{f}_{w})&=(d(z)-d(w))\textbf{f}_{w}+\sum_{s\in N(z)\setminus N[w]}\textbf{f}_{s}-\sum_{t\in X}\textbf{f}_{t}\\
 &\leq \lfloor\frac{k(k-1)}{4}\rfloor-k+2+\lfloor\frac{k(k-1)}{4}\rfloor+1\\
 &\leq \frac{k(k-1)}{2}-k+3.
\end{align*}
Since $d(z)=n-k$, we obtain
$$\textbf{f}_{z}-\textbf{f}_{w}\leq \frac{k^{2}-3k+6}{2(q(G)-n+k+1)}< \frac{k^{2}+6k+6}{2(q(G)-n+1)}.$$

\noindent\textbf{Case~2.} $Y_{1}\neq\emptyset$.

In this case, by Claim~3, we obtain that $\max_{v\in V(G)}\textbf{f}_{v}$ is attained by some vertex from $Y_{1}$, say, $z$. By the definition of $Y_{1}$, $z$ is adjacent to other vertices $G$.

If $\min_{v\in Y\cup Z}\textbf{f}_{v}$ is attained by some vertex from $Z_{2}$, say, $w$. Obviously, we have
$zw\in E(G)$ and $N(w)\setminus N[z]=\emptyset$. Furthermore, we observe that
$$d(z)-d(w)\leq ((n-1)-(n-k-\lfloor\frac{k(k-1)}{4}\rfloor-1))=\lfloor\frac{k(k-1)}{4}\rfloor+k,$$
and
$$|N(z)\setminus N[w]|\leq (k-1)+\lfloor\frac{k(k-1)}{4}\rfloor+1=\lfloor\frac{k(k-1)}{4}\rfloor+k.$$
Hence by \eqref{eq:11}, we have
\begin{align*}
(q(G)-d(z)+1)(\textbf{f}_{z}-\textbf{f}_{w})&=(d(z)-d(w))\textbf{f}_{w}+\sum_{s\in N(z)\setminus N[w]}\textbf{f}_{s}\\
 &\leq \lfloor\frac{k(k-1)}{4}\rfloor+k+\lfloor\frac{k(k-1)}{4}\rfloor+k\\
 &\leq \frac{k(k-1)}{2}+2k.
\end{align*}
Since $d(z)=n-1$, we obtain,
$$\textbf{f}_{z}-\textbf{f}_{w}\leq \frac{k^{2}+3k}{2(q(G)-n+1+1)}< \frac{k^{2}+6k+6}{2(q(G)-n+1)}.$$

If $\min_{v\in Y\cup Z}\textbf{f}_{v}$ is attained by some vertex from $Y_{2}$, say, $w$. Obviously, we have
$zw\in E(G)$ and $N(w)\setminus N[z]=\emptyset$. Furthermore, we observe that
$$d(z)-d(w)\leq ((n-1)-(n-1-\lfloor\frac{k(k-1)}{4}\rfloor-1))=\lfloor\frac{k(k-1)}{4}\rfloor+1,$$
and
$$|N(z)\setminus N[w]|\leq \lfloor\frac{k(k-1)}{4}\rfloor+1.$$
Hence by \eqref{eq:11}, we have
\begin{align*}
(q(G)-d(z)+1)(\textbf{f}_{z}-\textbf{f}_{w})&=(d(z)-d(w))\textbf{f}_{w}+\sum_{s\in N(z)\setminus N[w]}\textbf{f}_{s}\\
 &\leq \lfloor\frac{k(k-1)}{4}\rfloor+1+\lfloor\frac{k(k-1)}{4}\rfloor+1\\
 &\leq \frac{k(k-1)}{2}+2.
\end{align*}
So,
$$\textbf{f}_{z}-\textbf{f}_{w}\leq \frac{k^{2}-k+4}{2(q(G)-n+1+1)}< \frac{k^{2}+6k+6}{2(q(G)-n+1)}.$$

Claim~4 is proved.
\end{proof}

Hence we have
\begin{align}\label{eq:2}
\langle &\nonumber Q(G)\textbf{f},\textbf{f}\rangle-\langle Q((k-1)K_{1}+K_{n-k+1})\textbf{f},\textbf{f}\rangle\\
 &=\nonumber\sum_{y\in Y}d_{G[X\cup Y]}(y)\textbf{f}_{y}^{2}+\sum_{x\in X}d_{G[X\cup
 Y]}(x)\textbf{f}_{x}^{2}+2\sum_{\{x,y\}\in E(G[X\cup Y])}\textbf{f}_{x}\textbf{f}_{y}-\sum_{\{u,v\}\in
 E'}(\textbf{f}_{u}^{2}+\textbf{f}_{v}^{2}+2\textbf{f}_{u}\textbf{f}_{v})\\
 &=\nonumber\sum_{\{x,y\}\in E(G[X\cup Y])}(\textbf{f}_{x}+\textbf{f}_{y})^{2}-\sum_{\{u,v\}\in E'}(\textbf{f}_{u}+\textbf{f}_{v})^{2}\\
 &\nonumber\leq k(k-1)(1+\frac{k}{q(G)-k})^{2}-4|E'|(1-\frac{k^{2}+6k+6}{2(q(G)-n+1}))^{2}\\
 &< 0.
\end{align}

Here, the last inequality follows from $n\geq k^{4}+5k^{3}+2k^{2}+8k+12$ and $q(G)> 2n-2k-1$. We list the specific proof of \eqref{eq:2} in the last section. Note that $q((k-1)K_{1}+K_{n-k+1})=2n-2k> \langle Q(G)\textbf{f},\textbf{f}\rangle/\langle \textbf{f},\textbf{f}\rangle$. Hence,
$q(G)< 2n-2k$ for $G\in \mathcal{S}_{k}^{(2)}(n)$.

For $G\in \mathcal{T}_{k}^{(2)}(n)$, using the same definition of $Y_{1},Y_{2},Z_{1},Z_{2}$ as that of $\mathcal{S}_{k}^{(2)}(n)$, we can obtain the same conclusion as Claim~2 and Claim~3 by a similar method. Let $G\in \mathcal{T}_{k}^{(2)}(n)$ have the maximum signless Laplacian spectral radius. Let $\textbf{f}$ be the
eigenvector corresponding to $q(G)$. Furthermore, we assume $\max_{v\in V(G)}\textbf{f}_{v}=1$. Also, we
have
$$\max_{v\in Y\cup Z}\textbf{f}_{v}-\min_{u\in Y\cup Z}\textbf{f}_{u}\leq \frac{k^{2}+6k+6}{2(q(G)-n+1)}.$$
To prove this, we discuss the following two cases.

\noindent\textbf{Case~1.} $Y_{1}=\emptyset$.

In this case, $\max_{v\in V(G)}\textbf{f}_{v}$ is attained by some vertex from $Z_{1}$, say, $z$.

If $\min_{v\in Y\cup Z}\textbf{f}_{v}$ is attained by some vertex from $Z_{2}$, say, $w$.  Obviously, we have
$zw\in E(G)$, $N(w)\setminus N[z]=\emptyset$, $d(z)-d(w)\leq \lfloor\frac{k-1}{2}\rfloor+1$ and $|N(z)\setminus N[w]|\leq \lfloor\frac{k-1}{2}\rfloor+1$. Hence, by \eqref{eq:11}, we have
\begin{align*}
(q(G)-d(z)+1)(\textbf{f}_{z}-\textbf{f}_{w})&=(d(z)-d(w))\textbf{f}_{w}+\sum_{s\in N(z)\setminus N[w]}\textbf{f}_{s}\\
 &\leq k+1.
\end{align*}

If $\min_{v\in Y\cup Z}\textbf{f}_{v}$ is attained by some vertex from $Y_{2}$, say, $w$. Obviously, we have
$zw\in E(G)$, $N(w)\setminus N[z]=X$ and $|N(z)\setminus N[w]|\leq \lfloor\frac{k-1}{2}\rfloor+1$. Furthermore, we observe
$|d(z)-d(w)|\leq \lfloor\frac{k-1}{2}\rfloor+k+2.$
Hence by \eqref{eq:11}, we have
\begin{align*}
(q(G)-d(z)+1)(\textbf{f}_{z}-\textbf{f}_{w})&=(d(z)-d(w))\textbf{f}_{w}+\sum_{s\in N(z)\setminus N[w]}\textbf{f}_{s}-\sum_{t\in X}\textbf{f}_{t}\\
 &\leq \lfloor\frac{k-1}{2}\rfloor+k+2+\lfloor\frac{k-1}{2}\rfloor+1\\
 &\leq 2k+2.
\end{align*}

In both subcases, since $d(z)=n-k$, we have
$$\textbf{f}_{z}-\textbf{f}_{w}\leq \frac{k^{2}+6k+6}{2(q(G)-n+1)}.$$

\noindent\textbf{Case~2.} $Y_{1}\neq\emptyset$.

In this case, $\max_{v\in V(G)}\textbf{f}_{v}$ is attained by some vertex from $Y_{1}$. If $\min_{v\in Y\cup Z}\textbf{f}_{v}$ is attained by some vertex from $Z_{2}$, say, $w$. Obviously, we have
$zw\in E(G)$, $N(w)\setminus N[z]=\emptyset$. Furthermore, we observe that
$$d(z)-d(w)\leq ((n-1)-(n-k-\lfloor\frac{k-1}{2}\rfloor-1))=\lfloor\frac{k-1}{2}\rfloor+k,$$
and
$$|N(z)\setminus N[w]|\leq (k-1)+\lfloor\frac{k-1}{2}\rfloor+1=\lfloor\frac{k-1}{2}\rfloor+k.$$
Hence by \eqref{eq:11}, we have
\begin{align*}
(q(G)-d(z)+1)(\textbf{f}_{z}-\textbf{f}_{w})&=(d(z)-d(w))\textbf{f}_{w}+\sum_{s\in N(z)\setminus N[w]}\textbf{f}_{s}\\
 &\leq \lfloor\frac{k-1}{2}\rfloor+k+\lfloor\frac{k-1}{2}\rfloor+k\\
 &\leq 3k-1.
\end{align*}
Since $d(z)=n-1$, we obtain,
$$\textbf{f}_{z}-\textbf{f}_{w}\leq \frac{6k-2}{2(q(G)-n+1+1)}< \frac{k^{2}+6k+6}{2(q(G)-n+1)}.$$

If $\min_{v\in Y\cup Z}\textbf{f}_{v}$ is attained by some vertex from $Y_{2}$, say, $w$. Obviously, we have
$zw\in E(G)$ and $N(w)\setminus N[z]=\emptyset$. Furthermore, we observe that
$$d(z)-d(w)\leq ((n-1)-(n-1-\lfloor\frac{k-1}{2}\rfloor-1))=\lfloor\frac{k-1}{2}\rfloor+1,$$
and
$$|N(z)\setminus N[w]|\leq \lfloor\frac{k-1}{2}\rfloor+1.$$
Hence by \eqref{eq:11}, we have
\begin{align*}
(q(G)-d(z)+1)(\textbf{f}_{z}-\textbf{f}_{w})&=(d(z)-d(w))\textbf{f}_{w}+\sum_{s\in N(z)\setminus N[w]}\textbf{f}_{s}\\
 &\leq \lfloor\frac{k-1}{2}\rfloor+1+\lfloor\frac{k-1}{2}\rfloor+1\\
 &\leq k+1.
\end{align*}
So,
$$\textbf{f}_{z}-\textbf{f}_{w}\leq \frac{2k+2}{2(q(G)-n+1+1)}< \frac{k^{2}+6k+6}{2(q(G)-n+1)}.$$

Then for $G\in \mathcal{T}_{k}^{(2)}(n)$, repeating the argument for $G\in \mathcal{S}_{k}^{(2)}(n)$, we can
prove $q(G)< 2n-2k$ easily.

The proof is complete.
\hfill$\blacksquare$

Theorem \ref{th:6c0} can be stated as follows.

\noindent\begin{theorem}\label{th:6c4} Let $G$ be a graph of order $n\geq k^{4}+5k^{3}+2k^{2}+8k+12$ with minimum degree $\delta(G)\geq k$, where $k\geq 2$. If
$$q(G)\geq 2n-2k,$$
then $G$ is Hamilton-connected unless $G\in \mathcal{S}_{k}^{(1)}(n) \cup \mathcal{T}_{k}^{(1)}(n)$.
\end{theorem}

\noindent {\bf Proof.} By Theorem \ref{th:6cq}, we obtain
$$2n-2k\leq q(G)\leq \frac{2e(G)}{n-1}+n-2.$$
Note that $n\geq k^{4}+5k^{3}+2k^{2}+8k+12$, we have
\begin{align*}
e(G)&\geq \frac{(n-2k+2)(n-1)}{2}\\
 &= \frac{n^{2}-(2k+1)n+2n+2k-2}{2}\\
 &=\dbinom{n-k}{2}+\frac{2n-k^{2}+k-2}{2}\\
 &>\dbinom{n-k}{2}+k(k+1).
\end{align*}

By Theorem \ref{th:6c1}, $G$ is Hamilton-connected unless $G\subseteq S_{n}^{k}$ or $T_{n}^{k}$. Combining with
Lemmas \ref{th:6c2} and \ref{th:6c3}, we complete the proof.
\hfill$\blacksquare$

It is easy to see that if we do $k-2$ Kelmans operation on $T_{n}^{k}$, then we can obtain a proper
subgraph of $S_{n}^{k}$. Hence by the fact that signless Laplacian spectral radius decreases after deleting an edge in a connected graph, we obtain $q(S_{n}^{k})> q(T_{n}^{k})>q(K_{n-k+1})=2n-2k$. Then
we have the following corollary.

\noindent\begin{corollary}\label{co:6c0} Let $G$ be a graph of order $n\geq k^{4}+5k^{3}+2k^{2}+8k+12$ with minimum degree $\delta(G)\geq k$, where $k\geq 2$. If
$$q(G)\geq q(S_{n}^{k}),$$
then $G$ is Hamilton-connected unless $G= S_{n}^{k}$.
\end{corollary}

\section{Appendix}

\noindent\textbf{Proof of \eqref{eq:2}:}
Let $A=k(k-1)(1+\frac{k}{q(G)-k})^{2}$, $B=4|E'|(1-\frac{k^{2}+6k+6}{2(q(G)-n+1)})^{2}$. Since $q(G)\geq 2n-2k-1$,
we have
\begin{equation}\label{eq:3}
\max A=k(k-1)(1+\frac{k}{2n-3k-1})^{2}.
\end{equation}

If $k=4s$, then $|E'|=\lfloor \frac{4s(4s-1)}{4}\rfloor+1=s(4s-1)+1$. If $k=4s+1$, then $|E'|=\lfloor \frac{4s(4s+1)}{4}\rfloor+1=s(4s+1)+1$. In both cases, we have
$$\min B=(k(k-1)+4)(1-\frac{k^{2}+6k+6}{2(n-2k)})^{2}.$$
Combining with \eqref{eq:3}, we have
\begin{align*}
A&-B\leq \max A-\min B\\
&=k(k-1)(1+\frac{k}{2n-3k-1})^{2}-k(k-1)(1-\frac{k^{2}+6k+6}{2(n-2k)})^{2}-4(1-\frac{k^{2}+6k+6}{2(n-2k)})^{2}\\
 &=A_{1}+A_{2}+A_{3}-A_{4}-4,
\end{align*}
in which,
\begin{align*}
A_{1}&=\frac{2k^{3}-2k^{2}}{2n-3k-1};\\
A_{2}&=\frac{k^{4}-k^{3}}{4n^{2}-(12k+4)n+9k^{2}+6k+1};\\
A_{3}&=\frac{k^{4}+5k^{3}+4k^{2}+18k+24}{n-2k};\\
A_{4}&=\frac{k^{6}+11k^{5}+40k^{4}+72k^{3}+156k^{2}+252k+144}{4n^{2}-16kn+16k^{2}}.
\end{align*}
Since $n\geq k^{4}+5k^{3}+2k^{2}+8k+12$, we obtain that $A_{1}+A_{2}+A_{3}-A_{4}<4$. Hence, \eqref{eq:2} holds.

If $k=4s+2$, then $|E'|=\lfloor \frac{(4s+2)(4s+1)}{4}\rfloor+1= 4s^{2}+3s+1$. Hence,
$$B> 4(4s^{2}+3s+1)(1-\frac{k^{2}+3k+6}{n-2k})^{2}=(k(k-1)+2)(1-\frac{k^{2}+6k+6}{2(n-2k)})^{2}.$$
If $k=4s+3$, then $|E'|=\lfloor \frac{(4s+3)(4s+2)}{4}\rfloor+1= 4s^{2}+5s+2$. Hence,
$$B> 4(4s^{2}+5s+2)(1-\frac{k^{2}+3k+6}{n-2k})^{2}=(k(k-1)+2)(1-\frac{k^{2}+6k+6}{2(n-2k)})^{2}.$$
In both cases, combining with \eqref{eq:3}, we have
\begin{align*}
A&-B< \max A-B\\
&<k(k-1)(1+\frac{k}{2n-3k-1})^{2}-k(k-1)(1-\frac{k^{2}+6k+6}{2(n-2k)})^{2}-2(1-\frac{k^{2}+6k+6}{2(n-2k)})^{2}\\
 &=A'_{1}+A'_{2}+A'_{3}-A'_{4}-2,
\end{align*}
in which,
\begin{align*}
A'_{1}&=\frac{2k^{3}-2k^{2}}{2n-3k-1};\\
A'_{2}&=\frac{k^{4}-k^{3}}{4n^{2}-(12k+4)n+9k^{2}+6k+1};\\
A'_{3}&=\frac{k^{4}+5k^{3}+2k^{2}+6k+12}{n-2k};\\
A'_{4}&=\frac{k^{6}+11k^{5}+38k^{4}+48k^{3}+60k^{2}+108k+72}{4n^{2}-16kn+16k^{2}}.
\end{align*}
Since $n\geq k^{4}+5k^{3}+2k^{2}+8k+12$, we obtain that $A'_{1}+A'_{2}+A'_{3}-A'_{4}<2$. Hence, \eqref{eq:2} holds.

\hfill$\blacksquare$

\end{document}